\documentclass{amsart}

\input{epsf.tex}
\usepackage{amsfonts,amssymb,verbatim,amsmath,amsthm,latexsym,textcomp,amscd,hyperref,comment}
\usepackage{latexsym,amsfonts,amssymb,epsfig,verbatim}
\usepackage{amsmath,amsthm,amssymb,latexsym,graphics,textcomp}
\usepackage{graphicx}
\usepackage{color}
\usepackage{url}
\usepackage{stmaryrd}
\usepackage{tikz}
\usepackage{tikz-cd}
\usepackage{mathtools}
\usetikzlibrary{positioning}

\begin{document}

\newtheorem{theorem}{Theorem}[section]
\newtheorem{prop}[theorem]{Proposition}
\newtheorem{lemma}[theorem]{Lemma}
\newtheorem{cor}[theorem]{Corollary}
\newtheorem{prob}[theorem]{Problem}
\newtheorem{defn}[theorem]{Definition}
\newtheorem{notation}[theorem]{Notation}
\newtheorem{fact}[theorem]{Fact}
\newtheorem{conj}[theorem]{Conjecture}
\newtheorem{claim}[theorem]{Claim}
\newtheorem{example}[theorem]{Example}
\newtheorem{rem}[theorem]{Remark}
\newtheorem{assumption}[theorem]{Assumption}
\newtheorem{scholium}[theorem]{Scholium}

\newcommand{\map}{\rightarrow}
\newcommand{\C}{\mathcal C}
\newcommand\AAA{{\mathcal A}}
\def\AA{\mathcal A}

\def\L{{\mathcal L}}
\def\al{\alpha}
\def\A{{\mathcal A}}

\newcommand\GB{{\mathbb G}}
\newcommand\BB{{\mathcal B}}
\newcommand\DD{{\mathcal D}}
\newcommand\EE{{\mathcal E}}
\newcommand\FF{{\mathcal F}}
\newcommand\GG{{\mathcal G}}
\newcommand\HH{{\mathbb H}}
\newcommand\II{{\mathcal I}}
\newcommand\JJ{{\mathcal J}}
\newcommand\KK{{\mathcal K}}
\newcommand\LL{{\mathcal L}}
\newcommand\MM{{\mathcal M}}
\newcommand\NN{{\mathbb N}}
\newcommand\OO{{\mathcal O}}
\newcommand\PP{{\mathcal P}}
\newcommand\QQ{{\mathbb Q}}
\newcommand\RR{{\mathbb R}}
\newcommand\SSS{{\mathcal S}}
\newcommand\TT{{\mathcal T}}
\newcommand\UU{{\mathcal U}}
\newcommand\VV{{\mathcal V}}
\newcommand\WW{{\mathcal W}}
\newcommand\XX{{\mathcal X}}
\newcommand\YY{{\mathcal Y}}
\newcommand\ZZ{{\mathcal Z}}
\newcommand\hhat{\widehat}
\newcommand\flaring{{Corollary \ref{cor:super-weak flaring} }}
\def\Ga{\Gamma}
\def\Z{\mathbb Z}

\def\diam{\operatorname{diam}}
\def\dist{\operatorname{dist}}
\def\hull{\operatorname{Hull}}
\def\id{\operatorname{id}}
\def\Im{\operatorname{Im}}

\def\barycenter{\operatorname{center}}

\def\length{\operatorname{length}}
\newcommand\RED{\textcolor{red}}
\newcommand\BLUE{\textcolor{blue}}
\newcommand\GREEN{\textcolor{green}}
\def\mini{\scriptsize}

\def\acts{\curvearrowright}
\def\embed{\hookrightarrow}

\def\ga{\gamma}
\newcommand\la{\lambda}
\newcommand\eps{\epsilon}
\def\geo{\partial_{\infty}}
\def\bhb{\bigskip\hrule\bigskip}

\title[Acylindrical complexes of hyperbolic groups]{ Research announcement: A combination theorem for acylindrical complexes of hyperbolic groups and Cannon-Thurston maps}

\author{Pranab Sardar}
\address{Department of Mathematical Sciences,
	Indian Institute of Science Education and Research Mohali,
	Knowledge City, Sector 81, SAS Nagar,
	Punjab 140306,  India}

\email{psardar@iisermohali.ac.in}
\urladdr{https://sites.google.com/site/psardarmath/}

\author{Ravi Tomar}
\email{ravitomar547@gmail.com}
\address{Department of Mathematical Sciences,
	Indian Institute of Science Education and Research Mohali,
	Knowledge City, Sector 81, SAS Nagar,
	Punjab 140306,  India}

%\email{@iisermohali.ac.in}
%\urladdr{https://sites.google.com/site/psardarmath/}

%\thanks{MM is   supported by  the Department of Atomic Energy, Government of India, under project no.12-R\&D-TFR-14001.
%	MM is also supported in part by a Department of Science and Technology JC Bose Fellowship, CEFIPRA  project No. 5801-1, a SERB grant MTR/2017/000513, and %an endowment of the Infosys Foundation via the Chandrasekharan-Infosys Virtual Centre for Random Geometry.
%	This material is based upon work supported by the National Science Foundation
%	under Grant No. DMS-1928930 while MM participated in a program hosted
%	by the Mathematical Sciences Research Institute in Berkeley, California, during the
%	Fall 2020 semester. PS was partially supported by DST INSPIRE grant DST/INSPIRE/04/2014/002236 and DST MATRICS grant 
%MTR/2017/000485 of the Govt of India. }

\subjclass[2010]{20F65, 20F67 (Primary), 30F60(Secondary) }

\keywords{Complexes of groups, Cannon-Thurston map, hyperbolic groups, acylindrical action}

\date{\today}

\maketitle

%\tableofcontents
Complexes of groups (see \cite{bridson-haefliger} and \cite{haefliger}) are generalizations of graphs of groups \cite{serre-trees}.
Here we work with developable complexes of groups over finite simplicial complexes. Motivated by the celebrated combination 
theorem of Bestvina and Feighn \cite{BF} for graphs of hyperbolic groups, one asks if there is an analogous result for complexes of hyperbolic groups. One may start with a simplicial complex of groups $G(\YY)$ over a finite simplicial complex $Y$
 which is developable and such that 
(i) each face group is hyperbolic and (ii) for any two faces $\sigma\subset \tau$, the monomorphism $G_{\tau}\map G_{\sigma}$ 
is a qi embedding. However, under these assumptions it is in general unknown when the fundamental group will be hyperbolic
other than two extreme cases- (1) when all the monomorphisms $G_{\tau}\map G_{\sigma}$ have finite index image in $G_{\sigma}$ 
for any two faces $\sigma\subset \tau$ of $Y$ (\cite{mahan-sardar}) and (2) in case the action of the fundamental group of 
the complex of groups is acylindrical on the universal cover of $G(\YY)$ (\cite{martin1}). In both these cases additional 
hypotheses are needed. For instance, in \cite{martin1}, Martin works with simple complexes of groups, i.e. where there are
no twisting elements. However, using the result of \cite{martin-univ} one has the following theorem.

{\bf Theorem } \textup{(Martin)}\label{martin} {\em Let $G(\YY)$ be a developable complex of groups over a finite simplicial complex $Y$ such that the following holds:
	
	(1) Local groups are hyperbolic and local maps are quasi-isometric embeddings.
	
	(2) Universal cover of $G(\YY)$ is a CAT(0) hyperbolic space.
	
	(3) Action of $\pi_1(G(\YY))$ (fundamental group of $G(\YY))$ on universal cover is acylindrical.
	
Then $\pi_1(G(\YY))$ is hyperbolic and local groups are quasi-convex in $\pi_1(G(\YY))$.}

{\bf Remark:} In what follows we shall refer to the above theorem as `Martin's theorem'.

The idea behind the proof of the above theorem came from \cite{dahmani-comb}. Both Dahmani and Martin for the proofs of their theorems
constructed a candidate compact metric space on which there is a natural action of the
fundamental group of the complex of groups. Then they showed that the actions were uniform
convergence actions. Finally invoking Bowditch's theorem in \cite{bowditch-jams} they showed that the
fundamental group of the complex of groups is hyperbolic. In both cases for the proof to work
they needed all the face groups to be infinite. Contrary to this, in the first part of the paper we prove the following 
theorem in this connection.
\begin{theorem}
	Let $G(\YY)$ be a developable complex of groups over a finite simplicial complex $Y$ such that all the edge groups are finite. 
Suppose that the universal cover of $G(\YY)$ is a hyperbolic space. Then the fundamental group of $G(\YY)$ is hyperbolic relative to the
infinite vertex groups.
\end{theorem}
For the proof of the above theorem, we first show that the 1-skeleton of the universal cover is a fine graph; it is also hyperbolic by assumption. 
The theorem then follows from Bowditch's characterization (\cite{bowditch-relhyp}) of relatively hyperbolicity using fine hyperbolic graphs.

\begin{cor}
	Let $G(\YY)$ be a developable complex of groups such that the vertex groups are hyperbolic and the edge groups are all finite. 
Suppose the universal cover of $G(\YY)$ is a hyperbolic space. Then the fundamental group of $G(\YY)$, say $G$, is a hyperbolic group 
and vertex groups are quasiconvex in $G$.
\end{cor}
In the second part of the paper, motivated by the main theorem of \cite{mbdl2}, we ask the following question.

{\bf Question.} Let $G(\YY)$ be a complex of groups satisfying the hypotheses of Martin's theorem. Let $Y_1$ be a connected subcomplex of $Y$ and let $G(\YY_{1})$ be the subcomplex of groups obtained by restricting $G(\YY)$ to $Y_1$. Assume that the natural homomorphism from $H=\pi_1(G(\YY_{1}))$ to $G=\pi_1(G(\YY))$ is injective. Under what condition(s), is the group $H$ hyperbolic and when do we have a Cannon-Thurston map (see \cite{mahan-icm}) 
for the inclusion $H\rightarrow G$?

The following lemma was very useful for the proof of our next results. 

\begin{lemma} \textup{(\cite[Lemma 2.1]{mitra-trees})} 
Suppose $X,Y$ are  hyperbolic geodesic metric spaces and $f:Y\rightarrow X$ is an embedding. 
Then $f$ admits a Cannon-Thurston map if the following holds:
	
	Given $y_0\in Y$ there exists a non-negative function $M(N)$, such that $M(N)\rightarrow \infty $ as $N\rightarrow \infty$ and 
such that for all geodesic segments $\lambda$ lying outside $B(y_0,N)$ in $Y$, any geodesic segment in $X$ joining the end points of 
$f(\lambda)$ lies outside $B(f(y_0),M(N))$ in $X$. 
\end{lemma}
We note that the hypothesis of the above lemma implies that $f$ is a proper embedding. Also
Mitra's lemma gives a sufficient condition for the existence of a Cannon-Thurston map. It is not necessary (see \cite[Subsection 2.4]{mbdl2})
unless the spaces involved are proper metric spaces. However, we abstract out the hypothesis of the above lemma as a definition.

{\bf Definition.} {\em Suppose $X,Y$ are geodesic hyperbolic metric spaces and $f:Y\map X$ is any map.
We say that $f$ admits a Cannon-Thurston map in the strong sense if the following holds:

Given $y_0\in Y$ there exists a non-negative function $M(N)$, such that $M(N)\rightarrow \infty $ as $N\rightarrow \infty$ 
and such that for all geodesic segments $\lambda$ lying outside $B(y_0,N)$ in $Y$, any geodesic segment in $X$ joining end 
points of $f(\lambda)$ lies outside $B(f(y_0),M(N)$ in $X$.}

However, we have the following theorem.
\begin{theorem}\label{ct map}
	Let $G(\YY)$ be a complex of groups as in Martin's theorem and let $G(\YY_{1})$ be a subcomplex of groups as in the above question. 
Suppose that the natural homomorphism $H\rightarrow G$ is injective. Suppose $B_1$ is the universal cover of $G(\YY_{1})$ 
and $B$ is the universal cover of $G(\YY)$. Suppose $B_1$ is a CAT(0) hyperbolic space, and that there exists Cannon-Thurston 
map in the strong sense for the natural inclusion $B_1\rightarrow B$. Then the group $H$ is hyperbolic and there exists 
Cannon-Thurston map for the inclusion $H\rightarrow G$. Moreover, $H\rightarrow G$ is a qi embedding if and only if $B_1\rightarrow B$ admits an injective Cannon-Thurston map.
\end{theorem}
 
The idea of the proof of Theorem \ref{ct map} is the following. Hyperbolicity of the group $H$ is a direct consequence of Martin's theorem. 
Since both the complexes of groups are satisfying the hypotheses of Martin's theorem, we can explicitly construct the Gromov boundary for 
both $H$ and $G$ (see \cite{martin1}). Also from the construction of Gromov boundaries we have a natural $H$-equivariant map $\partial H\map \partial G$.
We prove that this map is continuous.\\

{\bf Some Examples.}

  Let $G(\YY)$ be a complex of groups over a polygon $Y$ satisfying all conditions of Martin's theorem.
In addition we assume that $Y$ has at least 4 edges. Let $Y_1$ be an edge of $Y$ and let $G(\YY_1)$ be 
the subcomplex of groups obtained by restricting $G(\YY)$ to $Y_1$. Let $H$ denote the fundamental group of 
$G(\YY_{1})$ and let $G$ denote the fundamental group of $G(\YY)$. Then we have:
 \begin{prop}
 	The natural homomorphism $H\map G$ is a qi embedding.
 \end{prop}
 The idea of the proof of the above proposition is the following. Let $T$ be the Bass-Serre tree of $G(\YY_1)$. 
We prove that the natural map from $T$ to the universal cover of $G(\YY)$ is an isometric embedding. Then the proposition
follows from the last part of Theorem \ref{ct map}.
 
The following special case of the proposition is worth mentioning.
 \begin{cor}
 	Let $G(\YY)$ be a complex of groups over a polygon $Y$ whose number of edges is $n\geq 4$. Assume that all the local groups are hyperbolic and all the local maps are qi embeddings and that the sum of group theoretical angles is less than $(n-2)\pi$. Also assume that the action of $\pi_1(G(\YY))$ on the universal cover of $G(\YY)$ is acylindrical. Let $Y_1$ be an edge in $Y$ and let $G(\YY_1)$ be the restriction of $G(\YY)$ to $Y_1$. Then $\pi_1(G(\YY))$ is a hyperbolic group and $\pi_1(G(\YY_1))$ is a quasi-convex subgroup of $\pi_1(G(\YY))$. 
 \end{cor}
For the definition of group theoretical angle at a corner of a polygon of groups, one is referred to \cite{gerstein-stalling}.


\begin{thebibliography}{Bow12}
	
	\bibitem[BF92]{BF}
	M.~Bestvina and M.~Feighn.
	\newblock A {C}ombination theorem for {N}egatively {C}urved {G}roups.
	\newblock {\em J. Diff. Geom., vol 35}, pages 85--101, 1992.
	
	\bibitem[BH99]{bridson-haefliger}
	M.~Bridson and A~Haefliger.
	\newblock Metric spaces of nonpositive curvature.
	\newblock {\em Grundlehren der mathematischen Wissenchaften, Vol 319,
		Springer-Verlag}, 1999.
	
	\bibitem[Bow98]{bowditch-jams}
	B.~H. Bowditch.
	\newblock A topological characterization of hyperbolic groups.
	\newblock {\em J. A. M. S. 11}, pages 643--667, 1998.
	
	\bibitem[Bow12]{bowditch-relhyp}
	B.~H. Bowditch.
	\newblock Relatively hyperbolic groups.
	\newblock {\em Internat. J. Algebra and Computation. 22, 1250016, 66pp}, 2012.
	
	\bibitem[Dah03]{dahmani-comb}
	Francois Dahmani.
	\newblock Combination of convergence groups.
	\newblock {\em Geometry and Topology, vol. 7}, pages 933--963, 2003.
	
	\bibitem[Hae91]{haefliger}
	Andr\'{e} Haefliger.
	\newblock Complexes of groups and orbihedra.
	\newblock In {\em Group theory from a geometrical viewpoint ({T}rieste, 1990)},
	pages 504--540. World Sci. Publ., River Edge, NJ, 1991.
	
	\bibitem[KS20]{mbdl2}
	S.~Krishna and P.~Sardar.
	\newblock {Pullbacks of metric bundles and Cannon-Thurston maps}.
	\newblock {\em https://arxiv.org/pdf/2007.13109.pdf}, 2020.
	
	\bibitem[Mar14]{martin1}
	A.~Martin.
	\newblock Non-positively curved complexes of groups and boundaries.
	\newblock {\em Geom. Topol. 18(1)}, pages 31--102, 2014.
	
	\bibitem[Mar15]{martin-univ}
	A.~Martin.
	\newblock Combination of universal spaces for proper actions.
	\newblock {\em J. Homotopy Relat. Struct. 10}, pages 803--820, 2015.
	
	\bibitem[Mit98]{mitra-trees}
	M.~Mitra.
	\newblock Cannon-{T}hurston {M}aps for {T}rees of {H}yperbolic {M}etric
	{S}paces.
	\newblock {\em Jour. Diff. Geom.48}, pages 135--164, 1998.
	
	\bibitem[Mj19]{mahan-icm}
	M.~Mj.
	\newblock {Cannon-Thurston maps}.
	\newblock {\em Proceedings of the International Congress of Mathematicians (ICM
		2018), World Scientific Publications, ISBN 978-981-3272-87-3}, pages
	885--917, 2019.
	
	\bibitem[MS12]{mahan-sardar}
	M.~Mj and P.~Sardar.
	\newblock {A Combination Theorem for metric bundles}.
	\newblock {\em Geom. Funct. Anal. 22, no. 6}, pages 1636--1707, 2012.
	
	\bibitem[Ser03]{serre-trees}
	J.-P. Serre.
	\newblock {\em Trees}.
	\newblock Springer Monographs in Mathematics. Springer-Verlag, Berlin, 2003.
	\newblock Translated from the French original by John Stillwell, Corrected 2nd
	printing of the 1980 English translation.
	
	\bibitem[Sta91]{gerstein-stalling}
	John~R. Stallings.
	\newblock Non-positively curved triangles of groups.
	\newblock In {\em Group theory from a geometrical viewpoint ({T}rieste, 1990)},
	pages 491--503. World Sci. Publ., River Edge, NJ, 1991.
	
\end{thebibliography}
\end{document}